\documentclass[a4paper,11pts]{article}
\usepackage{amsmath,amscd,amssymb,french}

\def \p {{\bf P^3}}
\def \pd {{\check \p}}
\def \I {{\bf I_3}}
\def \IT {{\bf I_3^1}}
\def \R {{\bf R}}
\def \ot {{\otimes}}
\def \V {{\check V}}
\def \F {{\bf F}}
\def \r {{\cal R}}
\def \fl {{\longrightarrow}}
\def \w {{\cal W}}
\def \prop {{\bf Proposition}}
\def \vs {\vskip}
\def \dm {{\bf D\'emonstration}}
\def \di {{\partial \I}}
\def \M {{\bf M}_{\p}(0,3,0)}
\def \Mn {{\bf M}_{\p}(0,n,0)}
\def \U {{\bf U}}
\def \du {{\partial \U}}

\begin{document}

{\centerline{\Large{\bf Eclatement de r\'eseaux de quadriques}}}
{\centerline{\Large{\bf et bord des instantons de degr\'e $3$}}}

\vs 0.2 cm

{\centerline{\large{\bf Blowing-up of nets of quadrics and boundary}}} 
{\centerline{\large{\bf of mathematical instantons of degree $3$}}}

\vs 0.3 cm

{\centerline{\large{\bf Nicolas PERRIN}}}
{\centerline{\large{\bf Ecole Normale Sup\'erieure}}}
{\centerline{\large{\bf 6 rue Einstein}}}
{\centerline{\large{\bf 92 160 Antony}}}
{\centerline{\large{\bf e-mail : \texttt{nperrin@clipper.ens.fr}}}}

\vs 0.5 cm

{\bf R\'esum\'e}:

\noindent
Dans leur article [1], L. Gruson et M. Skiti ont d\'ecrit une
application birationnelle de la vari\'et\'e $\I$ des instantons de
degr\'e $3$ vers les r\'eseaux de quadriques de $\pd$. Ils font ainsi
apparaitre deux composantes du bord de $\I$ associ\'ees l'une au
diviseur des r\'eseaux contenant une quadrique de $\pd$
d\'eg\'en\'er\'ee en deux plans et l'autre au diviseur des r\'eseaux
de L\" uroth. Dans cet
article, on d\'ecrit une composante irr\'eductible du bord de $\I$
comme le diviseur exceptionnel de l'\'eclatement du ferm\'e des
r\'eseaux de quadriques form\'e par les quadriques de rang $3$.

\vs 0.2 cm

{\bf Abstract} :

\noindent
In their article [1], L. Gruson and M. Skiti have constructed a
birationnal map from the variety $\I$ of mathematical instantons of
degree $3$ to the variety of nets of quadrics in $\pd$. They describe by
this way two irreducible componants of the boundary of $\I$ associated
to the divisor of nets which contain a two-plane degenerated quadric
and the divisor of L\" uroth nets. In this article
we describe an irreducible componante of the boundary of $\I$ as the
exceptionnal divisor of the blowing-up of the closed set of nets of
quadrics of rank $3$.

\vs 0.5 cm

{\bf Remerciements} :

\noindent
Je tiens \`a remercier ici mon directeur de th\`ese Laurent Gruson
pour toute l'aide qu'il m'a apport\'ee durant la pr\'eparation de ce travail.

\vs 0.5 cm

{\bf Introduction} :

Soit $V$ un espace vectoriel de dimension $4$ sur ${\bf C}$, on note
${\bf P}^3$ l'espace projectif ${\bf P}(V)$. Un instanton de degr\'e
$n$, est un \'el\'ement de $\Mn$ (l'espace des modules des faisceaux
semi-stables sans torsion de classes de chern $(0,n,0)$) qui est
localement libre, stable et qui v\'erifie la propri\'et\'e
cohomologique $h^1E(-2)=0$. On connait toutes les composantes
irr\'eductibles du bord des instantons de degr\'e $1$ et $2$ (voir par
exemple [3] pour le degr\'e $1$ et [2] pour le degr\'e $2$). On se
propose ici d'\'etudier une composante irr\'eductible du bord de la
vari\'et\'e $\I$ des instantons de degr\'e $3$.

L. Gruson et M. Skiti ont montr\'e dans [1] que la vari\'et\'e $\I$
des instantons de degr\'e $3$ est birationnelle \`a la vari\'et\'e
$\R={\rm{G}}(3,S^2\V)$ (qui param\'etrise les sous espace vectoriels
de dimension $3$ de $S^2\V$) des r\'eseaux de quadriques de
$\pd$. L'ouvert de $\I$ form\'e par les instantons sans droite
trisauteuse s'envoie birationnellement sur un ouvert $\R_4$ de $\R$
form\'e par les r\'eseaux $R$ de quadriques tels que l'application
$R\ot V\fl \V$ soit de rang $4$.

Cette description leur a permis d'identifier deux composantes
irr\'eductibles du bord de $\I$ donn\'ees par les hypersurfaces $\R'$
(respectivement $\R''$) des r\'eseaux contenant une quadrique
d\'ecompos\'ee en deux plans (respectivement des r\'eseaux de L\"
uroth). On va ici d\'ecrire une troisi\`eme composante du
bord. L'id\'ee consiste \`a \'eclater certaines sous vari\'et\'es de
$\R$ afin de pouvoir prolonger le morphisme des r\'eseaux vers les
faisceaux. On consid\`ere ainsi la sous vari\'et\'e $\F$ de la
Grassmannienne ${\rm{Grass}}(4,\r\ot V)$ des quotients de rang $4$ de
${\cal R}\otimes V$ au dessus de $\R$ (o\`u $\r$ est le sous
fibr\'e tautologique de $\R$) qui v\'erifie le fait que $\r\ot V\fl\w$
factorise l'application $\r\ot V\fl{\cal O}\otimes\V$ (on a ici not\'e $\w$ le
quotient tautologique de ${\rm{Grass}}(4,\r\ot V))$. Remarquons que $\pi
:\F \fl \R$ est un isomorphisme au dessus de $\R_4$. Notons $\R_i$ le
localement ferm\'e des r\'eseaux tels que la fl\`eche $\r\ot
V\stackrel{\psi}{\fl}\V$ est de rang $i$ et posons $\F_i=\pi^{-1}(\R_i)$.

\vskip 0.5 cm

{\bf Remarque} : Le sch\'ema $\F_3$ est irr\'eductible et r\'eduit. En
effet, on a un morphisme de $\F_3$ vers ${\rm{G}}(3,\V)$ qui a un
quotient $W$ de rang $4$ de ${\cal R}\otimes V$ associe l'image de la
compos\'ee ${\cal R}\otimes V\fl W\fl\V$ qui est de dimension $3$ car
on est dans $\F_3$. Mais alors la fibre de ce morphisme au dessus d'un
sous espace $K$ de dimension $3$ de $\V$ est donn\'ee par
${\rm{G}}(3,S^2K\oplus\V)$ donc on voit que $\F_3$ est donn\'e par
${\rm{G}}(3,S^2{\cal K}\oplus\V)$ o\`u ${\cal K}$ est le sous fibr\'e
tautologique de ${\rm{G}}(3,\V)$. Grace \`a cette description, on sait
que $\F_3$ est irr\'eductible et r\'eduit.

\vs 0.5 cm

\prop : \textit{Au voisinage de $\F_3$, le morphisme $\pi : \F\fl \R$ est
l'\'eclatement de $\R_3$.}

\vs 0.2 cm

\dm :
Prenons un r\'eseau $R$ de $\R_3$ et pla\c cons nous sur un ouvert affine
de $\R$ contenant ce point. On sait alors qu'il existe un mineur
$3\times 3$ de $R\ot V\fl\V$ qui est inversible. Ceci signifie que
l'on a un sous fibr\'e $Q$ de rang $3$ de ${\cal R}\ot V$ et un sous
fibr\'e $K$ de $\V$ tels que la restriction de $\psi$ \`a $Q$ est un
isomorphisme sur $K$. On peut alors se placer sur un ouvert affine de
$\R$ d'anneau $A$ tel que l'application $Q\ot A\stackrel{\psi}{\fl}
K\ot A$ est inversible en tout point de ${\rm{Spec}}(A)$. On a alors
le diagramme suivant de $A$-modules :
$$\begin{array}{ccc}
 Q & \stackrel{\sim}{\fl} & K \\
 \downarrow & & \downarrow \\
 \r\ot V & \stackrel{\psi}{\fl} & {\cal O}\otimes\V \\
 \downarrow & & \downarrow \\
 A^9 & \stackrel{\varphi}{\fl} & A \\
\end{array}$$
Le ferm\'e $\R_3$ dans l'ouvert affine ${\rm{Spec}}(A)$ est donn\'e
par le $0^{ieme}$ id\'eal de Fitting de $\varphi$ c'est \`a dire
l'annulation de $\varphi$. L'\'eclatement de ce lieu singulier est
alors un sch\'ema $X$ au dessus de ${\rm{Spec}}(A)$ solution du
probl\`eme universel suivant : si $X'$ est un sch\'ema muni d'un
morphisme $f$ vers ${\rm{Spec}}(A)$ tel que l'image de
$f^*\varphi : f^*(A^9){\fl} f^*(A)$ est un id\'eal inversible
alors le morphisme de $X$ vers ${\rm{Spec}}(A)$ se factorise par $f$.

Il reste \`a montrer que $\F$ au dessus de cet ouvert est \'egalement
solution de ce probl\`eme universel. Or au dessus de ${\rm{Spec}}(A)$,
le sch\'ema $\F$ v\'erifie la propri\'et\'e universelle suivante : Si
$X'$ est un sch\'ema muni d'un morphisme $f$ vers ${\rm{Spec}}(A)$ et
muni d'un module $W$ localement libre de rang $4$ tel que $f^*\psi :
f^*(\r\ot V)\fl f^*({\cal O}\otimes \V)$ se factorise par $W$, le
morphisme de $\F$ vers ${\rm{Spec}}(A)$ se factorise par $f$. Soit un
sch\'ema $X'$ muni d'un morphisme $f$ vers ${\rm{Spec}}(A)$ et tel que
l'image de $f^*\varphi : f^*(A^9)\fl f^*(A)$ est un id\'eal inversible
$I$, on construit alors un module localement libre de rang quatre $W$
tel que $f^*\psi$ se factorise par $W$. En effet, il suffit de poser
$W=f^*Q\oplus I$ qui est localement libre de rang $4$ et on a une
application $f^*(\r\ot V)\fl W$ qui factorise le morphisme
$f^*\psi$. La propri\'et\'e universelle de $\F$ au dessus de
${\rm{Spec}}(A)$ nous permet alors de dire que l'on a une application
de $\F$ vers $X'$ ce qui nous donne le r\'esultat.

\vs 0.5 cm

{\bf Remarque} : M. Skiti, dans un travail en pr\'eparation montre
qu'au voisinage de $\F_2$, le morphisme $\pi$ est l'\'eclatement de
$\R_2$. Il en d\'eduit une description de la sous vari\'et\'e (not\'ee
$\IT$) de $\I$ des instantons ayant une droite trisauteuse. L'application
birationnelle de $\R$ dans $\I$ s'\'etend \`a $\F$ de telle sorte qu'elle
devienne un isomorphisme sur un voisinage de $\F_2$ dans $\F$. Elle
identifie $\F_2$ et $\IT$. Ceci permet de d\'ecrire $\IT$ comme
l'\'eclatement de $\R_2$ dans $\R$. On va faire la m\^eme construction
avec $\F_3$.

\vs 0.5 cm

On \'etudie la sous vari\'et\'e $\di$ de $\M$ form\'ee par
les faisceaux sans torsion non localement libres dont le bidual est un
instanton de degr\'e $1$ et tels que le conoyau de l'injection
canonique dans ce bidual est une th\'eta-caract\'eristique (d\'ecal\'ee
de $2$) sur une conique lisse. Les faisceaux $E$ de cette famille sont donc
donn\'es par les noyaux de surjections $E''\fl\theta(2)$ o\`u $E''$
est un instanton de degr\'e $1$ et $\theta$ est une
th\'eta-caract\'eristique sur une conique lisse. Cette famille est
irr\'eductible de dimension $20$. Sur un ouvert de $\di$, les
faisceaux $E$ ont une cohomologie naturelle.

Sur l'ouvert $\U$ de $\M$ form\'e par les faisceaux \`a cohomologie
naturelle (donc minimale), on sait d\'efinir un morphisme vers
$\F$. Si $E$ est un tel faisceau, alors on a les \'egalit\'es
$h^1E(-1)=3$, $h^1E=4$ et $h^1E(1)=1$. Ainsi,
on d\'efinit un premier morphisme $f_0$ vers $\R$ qui prolonge celui de
[1] en associant \`a $E$ le r\'eseau $H^1E(-1)\fl H^1E(1)\ot
S^2\V$. Au dessus de ce r\'eseau, on associe \`a $E$ un quotient
de rang $4$ de $R\ot V$ donn\'e par $H^1E(-1)\ot V\fl H^1E$ qui nous
donne le morphisme $f$ souhait\'e. Ce morphisme est ainsi d\'efini sur
l'ouvert $\du$ de $\di$ form\'e par les faisceaux \`a cohomologie naturelle.

Il existe sur un ouvert de $\R$ (et donc sur un ouvert de
$\F$) un morphisme $g$ r\'eciproque (voir [1]). On va le prolonger
\`a un ouvert de $\F_3$.

\vs 0.5 cm

\prop : \textit{Le morphisme $f$ restreint \`a $\du$ est \`a valeurs dans
$\F_3$ et est dominant sur $\F_3$. De plus sur son image, on d\'efinit
une r\'eciproque $g$ \`a $f$.}

\vs 0.2 cm

\dm :
On commence par montrer que sur l'image de $\du$ on sait d\'efinir une
r\'eciproque \`a $f$. En effet, la r\'eciproque est donn\'ee de la
fa\c con suivante : soit $R$ le r\'eseau et $W$ le quotient de
$R\ot V$. Les applications lin\'eaires $R\ot V\to W$ et $W\to\V$ nous
donnent un complexe :
$$R\ot\Omega^2(2)\fl W\ot\Omega^1(1)\fl {\cal O}_{\p}$$
et la cohomologie au centre (d\'ecal\'ee de $1$) de ce complexe nous
donne le faisceau recherch\'e. En effet, si le r\'eseau $R$ et le quotient $W$ sont assez
g\'en\'eraux, le faisceau ainsi obtenu est dans $\M$. Soit $E$ dans
$\du$, la suite spectrale de Beilinson nous dit que le faisceau $E(1)$
est la cohomologie du complexe pr\'ec\'edent si on prend $R=H^1E(-1)$,
$W=H^1E$ et que l'on identifie $H^1E(1)$ \`a ${\bf C}$. Les applications
lin\'eaires sont les multiplications du module de Rao de $E$. Ainsi,
sur l'image de $\du$ par $f$ on a une r\'eciproque $g$ qui prolonge le
morphisme d\'efini sur $\R$.

Montrons maintenant que cette image est contenue dans $\F_3$. Pour
ceci, il suffit de voir que le r\'eseau associ\'e, qui est
$H^1E(-1)\fl H^1E(1)\ot S^2\V$ est tel que l'application
$H^1E(-1)\ot V\fl H^1E(1)\ot\V$ est de rang trois. Cette application se
d\'ecompose en deux applications : $H^1E(-1)\ot V\fl H^1E$ qui est
g\'en\'eriquement surjective  et $H^1E\fl
H^1E(1)\ot\V$. Mais si notre faisceau est donn\'e par la suite exacte 
$$0\fl E\fl E''\fl\theta(2)\fl 0$$
alors $H^1E$ s'identifie \`a $H^0\theta(2)$ et $H^1E(1)$ est un
quotient de rang $1$ de $H^0\theta(3)$. La multiplication du module de
Rao est donc nulle pour l'\'el\'ement $H$ de $V$ qui d\'efinit le plan
de la conique. Ainsi l'application $H^1E\to H^1E(1)\ot\V$ est de rang
trois et a pour image $(V/H{\check )}$. Par cons\'equent l'application
$H^1E(-1)\ot V\fl H^1E(1)\ot\V$ est aussi de rang trois car
$H^1E(-1)\ot V\fl H^1E$ est surjective.

On sait maintenant que $f$ a un morphisme r\'eciproque sur $f(\du)$
qui est contenu dans $\F_3$. On sait donc d\'ej\`a que l'image de
$\du$ est de dimension $20$. Mais alors, comme $\F_3$ est r\'eduit
irr\'eductible et que sa dimension est aussi $20$, on sait que $f(\du)$
contient un ouvert de $\F_3$ et $f\vert_{\du}$ est donc dominant sur $\F_3$.

\vs 0.5 cm

{\bf Corollaire} : \textit{Le prolongement de $g$ est birationnel
au voisinage de $\F_3$. Dans la description de $\I$ avec les
r\'eseaux, la vari\'et\'e $\di$ est donc le diviseur exceptionnel de
l'\'eclatement de $\R_3$ dans $\R$. La famille $\di$ forme une
composante irr\'eductible du bord de $\I$.}

\vs 0.2 cm

\dm :
On a vu que $g$ est d\'efini sur $f(\du)$ qui contient un ouvert de
$\F_3$. Ainsi, sur un voisinage (dans $\F$) de cet ouvert, le
morphisme $g$ est bien un morphisme r\'eciproque \`a $f$. Le morphisme
$g$ est bien birationnel au voisinage de $\F_3$.

De plus, on a vu que $\F\fl\R$ est l'\'eclatement de $\R_3$ au
voisinage de $\F_3$ donc $g$ identifie cette situation \`a celle de
$\I$ et $\di$. On voit donc que $\di$ est adh\'erente \`a $\I$ et que
dans la description de $\I$ avec les r\'eseaux, la vari\'et\'e $\di$
est le diviseur exceptionnel de l'\'eclatement de $\R_3$ dans $\R$ et
forme donc une composante irr\'eductible du bord de $\I$.

\vs 0.5 cm

{\bf Remarque} : On sait d\'ecrire les \'el\'ements de saut d'un
faisceau $E$ de $\di$. 

Les plans instables forment une courbe dans $\pd$ de degr\'e $6$ et de
genre $3$ ayant un point triple au point correspondant au plan de la
conique et dont la courbe des tris\'ecantes est trac\'ees sur le
complexe de droites associ\'e au bidual. 

Les droites bisauteuses forment une
courbe de degr\'e $8$ et de genre $3$ de la Grassmannienne qui est
r\'eunion d'une quintique rationnelle trac\'ee sur le complexe de
droites d\'efini par le bidual et d'une cubique trac\'ee dans le
$(\beta)$-plan des droites du plan de la conique. 

Il y a un lien entre ces deux courbes de saut. Les droites bisauteuses
sont les tris\'ecantes \`a la courbe des plans
instables. R\'eciproquement, la courbe des plans instables forme le
lieu triple de la surface r\'egl\'ee d\'ecrite par les droites
bisauteuses.

\vs 0.5 cm

{\bf R\'ef\'erences} :

\noindent
[1] L. Gruson et M. Skiti : $3$-instantons et r\'eseaux de
quadriques. Math. Ann. 298 (1994).

\vs 0.3 cm

\noindent
[2] M.S. Narashiman et G. Trautmann : Compactification of
$M_{\p}(0,2)$ and Poncelet pairs of conics. Pacific Journal of
Mathematics 145 (1990). 

\vs 0.3 cm

\noindent
[3] C. Okonek, M. Schneider et H. Spindler : Vector bunbles on
complex projective spaces. Basel, Boston, Stuttgart : Birkh\" auser 1980.

\end{document}